\documentclass[11pt]{article}
\usepackage{amsmath,amssymb,amsbsy,amsfonts,amsthm,latexsym,
            amsopn,amstext,amsxtra,euscript,amscd,amsthm}
\usepackage[cm]{fullpage}

\newtheorem{lem}{Lemma}
\newtheorem{lemma}[lem]{Lemma}

\newtheorem{thm}{Theorem}
\newtheorem{theorem}[thm]{Theorem}

\def\\{\cr}
\def\({\left(}
\def\){\right)}
\def\[{\left[}
\def\]{\right]}
\def\<{\langle}
\def\>{\rangle}

\begin{document}

\title{On ternary Egyptian fractions with prime denominator}
\author{Florian Luca\footnote{School of Mathematics, University of the Witwatersrand, Private Bag 3,Wits 2050 Johannesburg, South Africa, Research Group Algebraic Structures and Applications, King Abdulaziz University,
Jeddah, Saudi Arabia and Department of Mathematics, Faculty of Sciences, University of Ostrava, 30 Dubna 22, 701 03
Ostrava 1, Czech Republic} \& Francesco Pappalardi\footnote{Dipartimento di Matematica e Fisica, Universit\`a Roma Tre, Largo S. L. Murialdo 1, I--00141 Roma, Italy}}
\date{\today}

\pagenumbering{arabic}

\maketitle

\begin{abstract}
We prove upper and lower bound for the average  value  over primes $p$ of the number of positive integers $a$ such that the fraction $a/p$ can be written as the sum of three unit fractions.
\end{abstract}

An ``Egyptian fraction representation'' of a given rational $a/n$ is a solution in positive integers of the equation
$$\frac{a}{n}=\frac1{m_1}+\cdots\frac1{m_k}$$
In  case $k=2$  (resp. $k=3$) we  shall  say  it  is  a  \emph{binary (resp. ternary)  representation}.  A  variety  of  questions  about  these  representations  have  been  posed  and  studied.  Some  of these require them to be distinct but we shall not impose such a condition here. We refer  to  the  book  by Guy \cite{Guy} for a survey on this topic and an extensive list of references.

The object of our study is the following function: 
$$
A_k(n)=\#\left\{a\in {\mathbb N}: \frac{a}{n}=\frac{1}{m_1}+\frac{1}{m_2}+\cdots+\frac{1}{m_k},~m_1,m_2,\ldots,m_k\in {\mathbb N}\right\}.
$$
The case of binary Egyptian fractions was considered in \cite{pappa} where
it was shown that $A_2(n)\ll n^{o(1)}$ as $n\to\infty$ and that
$$
x\log^3x\ll\sum_{n\leq x}A_2(n)\ll x\log^3x.
$$
Note that some of the results in \cite{pappa} were improved in  \cite{HV}. 

The binary Egyptian fractions with prime denominators are significantly simpler. In fact it is quite easy to show that $A_2(p)=2 +\tau (p + 1)$
where $\tau$ is the divisor function. From this observation, it follows that
$$\sum_{p\le x}A_2(p)= \frac{315\zeta(3)}{2\pi^4}x+O\left(\frac{x}{\log x}\right)\qquad {\textrm{ as}}\qquad x\to\infty.$$
Here, we consider ternary Egyptian fractions and we study the average value of $A_3(p)$ as $p$ ranges over primes. In is shown in \cite{pappa} that $A_3(n)\ll n^{1/2+o(1)}$ as $n\to \infty$. We prove the following theorem.
\begin{theorem}\label{main}
We have 
$$
x(\log x)^3\ll \sum_{p\le x} A_3(p)\ll x(\log x)^5\quad\text{as}\ x\rightarrow\infty.
$$
\end{theorem}

As we shall see, the proof of the upper bound consists in estimating separately the contribution of fractions $m/p$ that admit a ternary Egyptian fraction expansion:
$$\frac mp=\frac1{m_1}+\frac1{m_2}+\frac1{m_3}.$$
with $p\mid\gcd(m_1,m_2)$ and $p\nmid m_3$ (Type I) and the contribution of those with $p\mid m_1$ and $p\nmid m_2m_3$  (Type II). 

The fraction of Type I are proven to contribute to $A_3(p)$ with $O(x\log^3x)$ and those of Type II with $O(x\log^5x)$.
We feel that, with some care, the $5$ in the latter log power should be lowered.
The proof of the above result was inspired by the paper of Elsholtz and Tao \cite{ElsTao} where, for a positive integer $n$, it is considered the number of integer solutions $x, y, z$ of the equation
$$\frac4n=\frac1x+\frac1y+\frac1z.$$

\section{Preliminaries}

We start with a description of $A_3(p)$. Recall the following result from \cite{BLP}. 

\begin{lemma}
\label{lem:main}
There is a representation of the reduced fraction $m/n$ (that is, $\gcd(m,n)=1$) as $m/n=1/m_1+1/m_2+1/m_3$ if and only if there are six positive integers $D_1,D_2,D_3,v_1,v_2,v_3$ with 
\begin{itemize}
\item[(i)] $\operatorname{lcm}[D_1,D_2,D_3]\mid n$, $\gcd(D_1,D_2,D_3)=1$;
\item[(ii)]   $v_1v_2v_3\mid D_1v_1+D_2v_2+D_3v_3$, and $\gcd(d_iv_i,v_j)=1$ for all $i\ne j$ with $\{i,j\}\in \{1,2\}$;
\item[(iii)] $m\mid (D_1v_1+D_2v_2+D_3v_3)/(v_1v_2v_3)$,
\end{itemize} 
and putting $E=\operatorname{lcm}[D_1,D_2,D_3]$, 
$f_1=n/E$, $f_2=(D_1v_1+D_2v_2+D_3v_3)/(mv_1v_2v_3)$ and $f=f_1f_2$, we have
\begin{equation}
\label{eq:x}
(m_1,m_2,m_3)=((E/D_1)v_2v_3f, (E/D_2)v_1v_3f, (E/D_3)v_1v_2f).
\end{equation}
\end{lemma}

Let us see the above lemma at work when $n=p$ is a prime. By condition (i) of the lemma, we first need $D_1,D_2,D_3$ such that $\operatorname{lcm}[D_1,D_2,D_3]\mid p$. This means that $\operatorname{lcm}[D_1,D_2,D_3]\in \{1,p\}$. 
The case in which $\operatorname{lcm}[D_1,D_2,D_3]=1$ leads to $D_1=D_2=D_3=1$ and now condition (ii)  shows that
$$
v_1v_2v_3\mid v_1+v_2+v_3.
$$
In particular, assuming $v_1\le v_2\le v_3$, we get $v_1v_2v_3\le 3v_3$, so $v_1v_2\le 3$. Thus, there are three possibilities for the pair $(v_1,v_2)$ and then since $v_3\mid v_1+v_2$, we infer that there are only the following three possibilities for $(D_1,D_2,D_3,v_1,v_2,v_3)$, namely $(1,1,1,1,1,1)$, $(1,1,1,1,1,2)$ and $(1,1,1,1,2,3)$.
Hence, $m\in\{1,2,3\}$.

Assume next that $\operatorname{lcm}[D_1,D_2,D_3]=p$. The situation $D_1=D_2=D_3=p$ is not possible since then the condition $(D_1,D_2,D_3)=1$ is not satisfied. 
Thus, not all $D_1,D_2,D_3$ are multiples of $p$. We then distinguish two cases. 

The first case is when there exists exactly one $D_i$ for $i\in \{1,2,3\}$ which equals  $p$. Say $D_1=p$. Then 
$$
v_1v_2v_3\mid pv_1+v_2+v_3\quad {\textrm{ and}}\quad m\mid (pv_1+v_2+v_3)/(v_1v_2v_3).
$$
In this case, $E=p$, so $E/D_1=1,~E/D_2=p$ and $E/D_3=p$. Furthermore, $f_1=1$ and $f_2=(pv_1+v_2+v_3)/(mp)$. Thus, $f=f_1f_2=f_2$. In addition, $v_1\mid v_2+v_3$ and
$$  
\frac{m}{p}=\frac{1}{v_2v_3 f}+\frac{1}{pv_1v_3f}+\frac{1}{pv_1v_2f}.
$$

In addition, $v_2$ and $v_3$ are coprime. 

Assume now that there are two indices $i,j\in \{1,2,3\}$ such that $D_i=D_j=p$. Say $D_2=D_3=p$. We then get that 
$$
v_1v_2v_3\mid v_1+pv_2+pv_3.
$$
Furthermore, $E=p$, $E/D_1=p,~E/D_2=1$, $E/D_3=1$, $f_1=1$. Moreover, $m\mid (pv_1+pv_2+v_3)/(v_1v_2v_3)$, $f_2=(pv_1+pv_2+v_3)/(mv_1v_2v_3)$, $f=f_2$ and 
$$
\frac{m}{p}=\frac{1}{pv_2v_3 f}+\frac{1}{v_1v_3 f}+\frac{1}{v_1v_2f}.
$$
If $p\mid v_1$, it follows that 
$$
\frac{m}{p}\le \frac{1}{p}+\frac{1}{p}+\frac{1}{p}\le \frac{3}{p},
$$
so again $m\in \{1,2,3\}$. Thus, suppose that $p\nmid v_1$. We then have from the fact that $v_1\mid p(v_2+v_3)$ together with the fact that $p\nmid v_1$, that $v_1\mid v_2+v_3$. In addition, $v_2$ and $v_3$ are coprime. 
 
To summarise, we proved the following lemma.
\begin{lemma}
If 
$$
\frac{m}{p}=\frac{1}{m_1}+\frac{1}{m_2}+\frac{1}{m_3}
$$
with positive integers $m_1,m_2,m_3$ and $\gcd(m,n)=1$, then either $m\in \{1,2,3\}$ or there exists positive integers $a,b$ with $\gcd(a,b)=1$, a positive integer $c$ with $c\mid a+b$ and a positive integer $u$ such that 
\begin{equation}
\label{eq:1}
\frac{m}{p}=\frac{1}{abu}+\frac{1}{pbcu}+\frac{1}{pacu}\qquad {\text{or}}\qquad \frac{m}{p}=\frac{1}{pabu}+\frac{1}{bcu}+\frac{1}{acu}.
\end{equation} 
\end{lemma}

We call the solutions from the left--hand side of \eqref{eq:1} \emph{solutions of Type I} and those from the right--hand side of equation \eqref{eq:1} \emph{solutions of Type II}.
The above lemma appears in many places (see \cite{ElsTao}, for example). However, we included the above proof of it since it can be deduced from the main result in \cite{BLP}.

The above lemma shows that either
\begin{equation}
\label{eq:1_1}
m=\frac{p+(a+b)/c}{abu}\qquad {\textrm{ or}}\qquad m=\frac{1+p(a+b)/c}{abu},
\end{equation}
where moreover $a$ and $b$ are coprime and $(a+b)/c$ is an integer. By symmetry, we always assume that $a\le b$. 

 Recall that the goal in order to estimate 
$$
\sum_{p\le x} A_3(p).
$$
That is, to count pairs $(m,p)$ with $p\le x$  such that $m/p$ can be written as a Egyptian fraction with three summands. If $\gcd(m,p)>1$, then $p\mid m$. Thus, 
$m/p=k$ is an integer and since it equals $1/m_1+1/m_2+1/m_3$ for some positive integers $m_1,m_2,m_3$, we have that $k\in \{1,2,3\}$ and then $m\in \{p,2p,3p\}$. 
Thus, 
$$
\sum_{p\le x} A_3(p)=\sum_{p\le x} A_3^*(p)+O(\pi(x))=\sum_{p\le x} A_3^*(x)+O(x/\log x).
$$
Thus, it suffices to count pairs $m/p$ with $\gcd(m,p)=1$. We start with lower bounds.

\section{Lower bound}

To prove the lower bound we count fractions $m/p$ with $m\not\in \{1,2,3,p,2p,3p\}$ arising from solutions of Type I for a large $x$ with the following property:
\begin{itemize}
\item[(i)] $a\in [x^{1/200}, x^{1/100}]$, $\tau(a)<(\log x)^4$;
\item[(ii)] $b\in [x^{1/20},x^{1/10}]$, $\tau(b)<(\log x)^4$;
\item[(iii)] $\tau(a+b)<(\log x)^4$;
\item[(iv)] $c\mid a+b$, $c\in [x^{1/200},x^{1/100}]$;
\item[(v)] $u\in [x^{1/200},x^{1/100}]$ is coprime to $a+b$ and $\tau(u)<(\log x)^4$.
\end{itemize}
We let ${\mathcal A}(x)$ be the set of quadruples $(a,b,c,u)$ with the above property. For such a quadruple $(a,b,c,u)\in {\mathcal A}(x)$, we have 
$$
m=\frac{p+(a+b)/c}{abu}.
$$
Thus, $p\equiv d^*\pmod {abu}$, where $d^*:=d(a,b,c)$ is the residue class  of the number $-(a+b)/c$ modulo $abu$. Note that $(a+b)/c$ is coprime to $ab$ because $a$ and $b$ are coprime and $(a+b)/c$ is coprime to $u$ by construction. 
Before we dig into getting a lower bound, we ask whether distinct quadruples $(a,b,c,u)$ as above give rise to distinct fractions $m/p$. Well, let us suppose that they do not and that there are 
$(a,b,c,u)\ne (a_1,b_1,c_1,u_1)$ such that $m/p=m_1/p_1$. Since $m\not\in \{p,2p,3p\}$ it follows that $m/p$ is not an integer. Hence, $m_1/p_1$ is not an integer either, so $m/p=m_1/p_1$ entails $p=p_1$ and  $m=m_1$. 
So, we get
$$
\frac{p+(a+b)/c}{abu}=\frac{p+(a_1+b_1)/c_1}{a_1b_1u_1}.
$$
In turn this gives
$$
p(abu-a_1b_1u_1)=((a+b)/c)a_1b_1u-((a_1+b_1)/c_1)abu.
$$
Assume first that the right--hand side above is nonzero. The left--hand side is nonzero also. Then $p$ is a divisor of $|((a+b)/c)a_1b_1u_1-((a_1+b_1)/c_1)abu|$, a nonzero number of size at most 
$x^{O(1)}$, which therefore has at most $O(\log x)$ prime factors. Further, the eight-tuple $(a,b,c,u,a_1,b_1,c_1,u_1)$ can be chosen in at most 
$$
x^{2(1/10+1/100+1/100)+o(1)}<x^{1/4}\quad {\textrm{ ways for large}}\quad x.
$$ 
Thus, there are at most $x^{1/4}\log x$ 
primes $p$ that can appear in that way, and for each such prime $p$ we have $A_3(p)\le p^{1/2+o(1)}$ as $p$ tends to infinity by one of the results from \cite{pappa}. Thus, for large $x$, there are at most $x^{1/4+1/2+o(1)}<x^{4/5}$ pairs $(m,p)$ for large $x$ with the property that
$m/p$ arises from two different quadruples $(a,b,c,u)$ and $(a_1b_1,c_1,u_1)$ as above for which $abu\ne a_1b_1u_1$. Since we are shooting for a lower bound of $\gg x(\log x)^3$, these pairs are negligible for the rest of the argument. 

Assume next that $abu=a_1b_1u_1$. In this case, we also get
\begin{equation}
\label{eq:2}
\frac{a+b}{c}=\frac{a_1+b_1}{c_1}.
\end{equation}
Since $abu=a_1b_1u_1$ and $\gcd(b,a_1u_1)\le a_1u_1\le x^{2/100}=x^{1/50}$, it follows that 
$$
\gcd(b,b_1)\ge bx^{-1/50}\ge x^{1/20-1/50}=x^{3/100}.
$$
Reducing equation \eqref{eq:2} modulo $\gcd(b,b_1)$, we get that
$$
ac_1-a_1c\equiv 0\pmod {\gcd(b,b_1)}.
$$
The left--hand side is an integer in absolute value at most 
$$
|ac_1-a_1c|\le \max\{ac_1,a_1c\}\le x^{1/100+1/100}=x^{2/100}<x^{3/100}\le \gcd(b,b_1).
$$
It thus follows that $ac_1=a_1c$ so $a/c=a_1/c_1$. Since $\gcd(a,c)=\gcd(a_1,c_1)=1$ (because $c\mid a+b$ and $c_1\mid a_1+b_1$ and $\gcd(a,b)=\gcd(a_1,b_1)=1$), it follows that $a=a_1,~c=c_1$. Now \eqref{eq:2}
implies $b=b_1$ and the equality $abu=a_1b_1u_1$ implies now that $u=u_1$, so $(a,b,c,u)=(a_1,b_1,c_1,u_1)$, a contradiction. The above argument shows that
$$
\sum_{p\le x} A_3(p)\ge \sum_{(a,b,c,u)\in {\mathcal A}(x)} \sum_{\substack{p\le x\\ p\equiv d^*\bmod {abu}}} 1+O(x^{4/5}),
$$
and it remains to deal with the first sum which equals:
$$\sum_{(a,b,c,u)\in {\mathcal A}(x)} \pi(x,abu,d^*).$$
For this, we use the Bombieri-Vinogradov Theorem. Note that we are counting primes $p$ in a certain arithmetic progression of ratio $abu<x^{1/10+1/100+1/100}=x^{3/25}<x^{1/3}$. The Bombieri-Vinogradov Theorem tells us that for every $A$, we have 
\begin{equation}
\label{eq:3}
\sum_{Q\le x^{1/3}} \max_{\substack{y\le x\\ 1\le d\le Q, (d,Q)=1}} \left|\pi(x,Q,d)-\frac{\pi(x)}{\varphi(Q)}\right|\ll_A\frac{x}{(\log x)^A}.
\end{equation}
For us, we will take $Q=abu$. However, given $Q$, there are many ways to choose $(a,b,u)$ and then even more ways to choose $c$. Well, let us count how many ways there are. 
We have $\tau(Q)=\tau(abu)\le \tau(a)\tau(b)\tau(u)\le (\log x)^{12}$ by properties (i), (ii), (v). Thus, the triple $(a,b,u)$ with $abu=Q$ and $\gcd(a,b)=\gcd(a+b,u)=1$ can be chosen 
in at most $\tau(Q)^2\le (\log x)^{24}$ ways. Having chosen $(a,b,u)$, we have that 
$c\mid a+b$, so $c$ can be chosen in at most $\tau(a+b)\le (\log x)^4$ ways. Hence, $(a,b,c,u)$ can be chosen in at most $(\log x)^{28}$ ways. Note that $(a,b,c)$ determine $d^*$ uniquely via $d^*\equiv -(a+b)/c\pmod Q$. 
Thus, for each $Q=abu$, we have at most $(\log x)^{28}$ values of $d^*$. Taking $A=30$ in \eqref{eq:3}, we get that 
\begin{equation}
\label{eq:50}
 \sum_{(a,b,c,u)\in {\mathcal A}(x)}  \pi(x,abu,d^*)=\pi(x)\sum_{(a,b,c,u)\in {\mathcal A}(x)} \frac{1}{\varphi(abu)}+O\left(\frac{x}{(\log x)^2}\right).
 \end{equation}
 We need to deal with the sum on the right--hand side above. Putting, ${\mathcal I}:=[x^{1/200},x^{1/100}]$, ${\mathcal J}:=[x^{1/20},x^{1/10}]$, 
 $\tau_{\mathcal I}(n):=\sum_{\substack{d\mid n\\ d\in {\mathcal I}}} 1$, we have
 \begin{equation}
 \label{eq:4}
 \sum_{(a,b,c,u)\in {\mathcal A}(x)} \frac{1}{\varphi(abu)}=\sum_{\substack{a\in {\mathcal J}\\ \tau(a)<(\log x)^4}}\sum_{\substack{b\in {\mathcal I}\\ (a,b)=1\\ \max\{\tau(b),\tau(a+b)\}<(\log x)^4}}
 \frac{\tau_{\mathcal I}(a+b)}{\varphi(a)\varphi(b)} \sum_{\substack{u\in {\mathcal I}\\ (u,a+b)=1\\ \tau(u)<(\log x)^4}} \frac{1}{\varphi(u)}.
 \end{equation}
 It is easy to sum up reciprocals. What gets in the way are the extra conditions, which are coprimality and the restriction on the size of the divisors functions. Let us start with the inner sum. Since $\varphi(u)\le u$, we have
 \begin{eqnarray}
 \label{eq:5}
 \sum_{\substack{u\in {\mathcal I}\\ (u,a+b)=1\\ \tau(u)<(\log x)^4}} \frac{1}{\varphi(u)} & \ge & \sum_{\substack{u\in {\mathcal I}\\ (u,a+b)=1\\ \tau(u)<(\log x)^4}} \frac{1}{u}\nonumber\\
 & \ge & \sum_{\substack{u\in {\mathcal I}\\ (u,n)=1}} \frac{1}{u}-\sum_{\substack{u\in {\mathcal I}\\ \tau(u)\ge (\log x)^4}} \frac{1}{u}\nonumber\\
 & := & S_{1,n}-S_2
 \end{eqnarray}
 with $n:=a+b$. We need an upper bound on $S_2$ and a lower bound on $S_{1,n}$. We start with the upper bound on $S_2$. Since 
 $$
 \sum_{u\le t} \tau(u)\le t\log t\ll t\log x\quad {\textrm{ for~all}}\quad t\in [x^{1/200},x^{1/100}],
 $$
 it follows that if we set ${\mathcal D}:=\{u: \tau(u)\ge (\log x)^4\}$ and put $D(t):={\mathcal D}\cap [1,t]$, then 
 \begin{equation}
 \label{eq:99}
 \#{\mathcal D}(t)\ll t/(\log x)^3\quad {\textrm{ uniformly for}}\quad t\in {\mathcal I}.
 \end{equation}
 Thus, by the Abel summation formula, 
 \begin{eqnarray}
 \label{eq:100}
S_2 & = &  \sum_{u\in {\mathcal D}\cap {\mathcal I}}\frac{1}{u}\nonumber\\
& = & \frac{\#{\mathcal D}(x^{1/100})}{x^{1/100}}-\frac{\#{\mathcal D}(x^{1/200})}{x^{1/200}}-\int_{x^{1/200}}^{x^{1/100}} \#{\mathcal D}(t) \left(-\frac{1}{t^2}\right) dt\nonumber\\
 & \ll & \int_{x^{1/200}}^{x^{1/100}} \frac{1}{t(\log x)^3} dt+\frac{1}{(\log x)^3}\nonumber\\
 & \ll & \frac{1}{(\log x)^2}.
 \end{eqnarray}
 We now discuss $S_{1,n}$. Clearly, if we take ${\mathcal D}=\{u: (n,u)=1\}$,  and put ${\mathcal D}(t)={\mathcal D}\cap [1,t]$, we have
 $$
 \#{\mathcal D}(t)=\sum_{d\mid n} \mu(d)\left\lfloor \frac{t}{d}\right\rfloor=\sum_{d\mid n} \mu(d)\left(\frac{t}{d}+O(1)\right)=\frac{\varphi(n)}{n} t+O(\tau(n)).
 $$
 In particular, 
 \begin{eqnarray*}
 S_{1,n} & = & \sum_{\substack{u\in {\mathcal I}\\ u\in {\mathcal D}}} \frac{1}{u}\\
 & = & \frac{\#{\mathcal D}(x^{1/100})}{x^{1/100}}-\frac{\#{\mathcal D}(x^{1/200})}{x^{1/200}}-\int_{x^{1/200}}^{x^{1/100}} \#{\mathcal D}(t) \left(-\frac{1}{t^2}\right) dt\\
 & = & \int_{x^{1/200}}^{x^{1/100}} \left(\frac{\varphi(n)}{n} t+O(\tau(n))\right)\frac{1}{t^2} dt+O\left(\frac{\tau(n)}{x^{1/200}}\right)\\
 & = & \frac{\varphi(n)}{n}\int_{x^{1/200}}^{x^{1/100}} \frac{dt}{t}+O\left(\tau(n)\int_{x^{1/200}}^{x^{1/100}} \frac{dt}{t^2}+\frac{\tau(n)}{x^{1/200}}\right)\\
 & \gg & \frac{\varphi(n)}{n}\log x+O(\tau(n)x^{-1/200})\\
 & \gg & \frac{\varphi(n)}{n}\log x+O(x^{-1/201}),
 \end{eqnarray*}
 where we use the fact that $\tau(n)=x^{o(1)}$ for $x\to\infty$, and in particular $\tau(n)x^{-1/200}\ll x^{-1/201}$. Since $(\varphi(n)/n)\log x\gg \log x/\log\log x$ and $x^{-1/201}=o(1)=o(\log x/\log\log x)$
as $x\to\infty$, it follows that in the above estimate, we may neglect the second term in the right--most side. Hence, 
$$
S_{1,n}\gg \frac{\varphi(n)\log x}{n}.
$$
We thus get that
\begin{equation}
\label{eq:6}
S_{1,n}-S_2\gg \frac{\varphi(n)}{n}\log x-O\left(\frac{1}{(\log x)^2}\right)\gg \frac{\varphi(n)}{n}\log x.
\end{equation}
Thus, using \eqref{eq:6} into \eqref{eq:5}, \eqref{eq:4} becomes 
\begin{equation}
\label{eq:60}
 \sum_{(a,b,c,u)\in {\mathcal A}(x)} \frac{1}{\varphi(abu)} \gg \log x \sum_{\substack{a\in {\mathcal J}\\ \tau(a)<(\log x)^4}}\sum_{\substack{b\in {\mathcal I}\\ (a,b)=1\\ \max\{\tau(b),\tau(a+b)\}<(\log x)^4}}
 \frac{\tau_{\mathcal I}(a+b)\varphi(a+b)}{ab(a+b)}.
 \end{equation}
 Now observe that $n:=a+b$ is in the interval $[x^{1/20}+x^{1/200},x^{1/10}+x^{1/100}]$. We shrink this to ${\mathcal J}_1:=[x^{1/19},x^{1/10}]$ and consider $n=a+b\in {\mathcal J}_1$ with $a\in [x^{1/200},x^{1/100}]$ coprime to $n$. In fact, $a$ is coprime to $n$ if and only if $a$ coprime to $b$. Further, $b=n-a>n/2$. So, the sums in the right--hand side 
 of \eqref{eq:60} above 
 exceed
 $$
\gg  \sum_{\substack{n\in {\mathcal J}_1\\ \tau(n)<(\log x)^4}} \frac{\tau_{\mathcal I}(n)\varphi(n)}{n^2}\sum_{\substack{a\in {\mathcal I}\\ (a,n)=1\\ \max\{\tau(a),\tau(n-a)\}<(\log x)^4}}\frac{1}{a}.
 $$
The extra condition $\tau(n-a)<(\log x)^4$ is a translation of the condition $\tau(b)<(\log x)^4$ with the new notations. We get that for fixed $n$, the inner sum satisfies
\begin{eqnarray*}
 \sum_{\substack{a\in {\mathcal I}\\ (a,n)=1\\ \max\{\tau(a),\tau(n-a)\}<(\log x)^4}}\frac{1}{a} & \ge & \sum_{\substack{a\in {\mathcal I}\\ (a,n)=1}}\frac{1}{a}-\sum_{\substack{a\in {\mathcal I}\\ \tau(a)>(\log x)^4}}\frac{1}{a}-
 \sum_{\substack{a\in {\mathcal I}\\ \tau(n-a)>(\log x)^4}}\frac{1}{a}\\
 & : = & S_{1,n}-S_2-S_3,
 \end{eqnarray*}
say. By the previous arguments, we have that $S_{1,n}\gg (\varphi(n)/n)\log x$, and $S_2\ll (\log x)^{-2}$. It remains to deal with $S_3$. Luckily, this has been done in \cite{ElsTao}. Namely, Corollary 7.4 in \cite{ElsTao}, shows that
 uniformly for $t\in {\mathcal I}$, we have that 
 $$
 \sum_{a\le t} \tau(n-a)\ll t\log t\ll t\log x.
 $$
 Thus, putting ${\mathcal D}:=\{a: \tau(n-a)\ge (\log x)^4\}$, and ${\mathcal D}(t):={\mathcal D}\cap [1,t]$, we have that
 $$
 \#{\mathcal D}(t) \ll t/(\log x)^3\quad {\textrm{ uniformly~for}}\quad t\in {\mathcal I}.
 $$
 This is enough, via the Abel summation formula as in the argument used to derive \eqref{eq:100} from \eqref{eq:99}, to deduce that 
 $$
 S_3\ll  (\log x)^{-2}.
 $$
 Hence, we get that $S_{1,n}-S_2-S_3\gg (\varphi(n)/n)\log x$, so that
 $$
 \sum_{(a,b,c,u)\in {\mathcal A}(x)} \frac{1}{\varphi(abu)}\gg (\log x)^2\sum_{\substack{n\in {\mathcal J}_1\\ \tau(n)<(\log x)^4}}\frac{\tau_{\mathcal I}(n)\varphi(n)^2}{n^3}.
 $$
 Lastly we need to worry about numbers with a bounded number of divisors, so we write the last sum as
 \begin{eqnarray*}
 \sum_{\substack{n\in {\mathcal J}_1\\ \tau(n)<(\log x)^4}}\frac{\tau_{\mathcal I}(n)\varphi(n)^2}{n^3} & \ge & \sum_{n\in {\mathcal J}_1}\frac{\tau_{\mathcal I}(n)\varphi(n)^2}{n^3}-\sum_{\substack{n\in {\mathcal J}_1\\ \tau(n)\ge (\log x)^4}}
 \frac{\tau_{\mathcal I}(n)}{n}\\
 & : = & S_1-S_2.
 \end{eqnarray*}
To bound $S_2$ we note that by writing $n=dv$ for some divisor $d\in {\mathcal I}$, by changing the order of summation, we have
\begin{eqnarray*}
S_2 & = &\sum_{\substack{n\in {\mathcal J}_1\\ \tau(n)\ge (\log x)^4}}
 \frac{\tau_{\mathcal I}(n)}{n}\\
 &  = & \sum_{d\in {\mathcal I}} \sum_{\substack{n\in {\mathcal J}_1\\ d\mid n\\ \tau(n)\ge (\log x)^4}} \frac{1}{n}\\
& = & \sum_{d\in {\mathcal I}}\sum'_{x^{1/19}/d\le v\le x^{1/10}/d} \frac{1}{dv},
\end{eqnarray*}
where the prime $'$ notation encodes the condition that $\tau(dv)\ge (\log x)^4$. Since $\tau(dv)\le \tau(d)\tau(v)$, it follows that either $\tau(d)\ge (\log x)^2$ or $\tau(v)\ge (\log x)^2$. Retaining this condition for either $d$ or $v$ and summing up trivially over the other parameter, we get that 
\begin{equation}
\label{eq:1000}
S_2\ll \log x\sum_{\substack{x^{1/200}\le d\le x^{1/10}\\ \tau(d)\ge (\log x)^2}}\frac{1}{d}.
\end{equation}
The counting function of the last set ${\mathcal D}:=\{d: \tau(d)\ge (\log x)^2\}$ satisfies the inequality
$$
\#{\mathcal D}(t)\ll t/\log x\quad {\textrm{ uniformly~in}}\quad t\in [x^{1/200},x^{1/20}],
$$
where as usual ${\mathcal D}(t)={\mathcal D}\cap [1,t]$. By the Abel summation formula, we get that
$$
\sum_{\substack{x^{1/200}\le d\le x^{1/10}\\ \tau(d)\ge (\log x)^2}}\frac{1}{d}=O(1),
$$
showing via \eqref{eq:1000} that $S_2=O(\log x)$. Finally, 
\begin{eqnarray*}
S_1 & = & \sum_{n\in {\mathcal J}_1} \frac{\tau_{\mathcal I}(n)\varphi(n)^2}{n^3} \\
& = & \sum_{d\in {\mathcal I}} \sum_{\substack{n\in {\mathcal J}_1\\ d\mid n}} \frac{\varphi(n)^2}{n^3}\\
& = & \sum_{d\in {\mathcal I}} \sum_{x^{1/19}/d\le v\le x^{1/10}/v} \frac{\varphi(dv)^2}{(dv)^3}\\
& \gg & \left(\sum_{x^{1/200}\le d\le x^{1/100}} \frac{\varphi(d)^2}{d^3}\right)\left( \sum_{x^{1/19-1/200}\le v\le x^{1/10-1/100}} \frac{\varphi(v)^2}{v^3}\right)\\
& \gg & (\log x)^2.
\end{eqnarray*}
This shows that $S_1-S_2\gg (\log x)^2$, and therefore that 
$$
\sum_{p\le x} A_3(p)\gg \pi(x)(\log x)^4+O\left(\frac{x}{(\log x)^2}+x^{4/5}\right)\gg x(\log x)^3.
$$

\section{Upper bound}

We shall bound the sum in the statement of Theorem~\ref{main} restricted to primes $p$ that admit solutions of Type I and Type II separately. 

\subsection{Type I solutions} In this case, from (\ref{eq:1_1}), we have
\begin{equation}
\label{eq:****}
m=\frac{p+(a+b)/c}{abu}.
\end{equation}
By symmetry, we may assume that $a\le b$. We may also assume that $m\ge (\log x)^4$, otherwise there are only $O(\pi(x)(\log x)^4)=O(x(\log x)^3)$ pairs  of positive integers $(m,p)$ with $p\le x$ and $m\le (\log x)^4$, and this bound is acceptable for us. 
Thus, $abu\ll x/(\log x)^4$. Let $\delta>0$ to be fixed later. 

\medskip

{\bf Case 1.} \textit{Assume that $abu\le x^{1-\delta}$}.

\medskip

Let $f_1(p)$ be the number of $m$ arising in this way from some $p$. Then fixing $abu$ and $c\mid a+b$, we need to count the number of primes $p\le x$ with $p\equiv d^*\pmod {abu}$, where $d^*$ is the congruence class of $-(a+b)/c$ modulo $abu$. Clearly, $(a+b)/c$ and $ab$ are coprime. The event that $u$ is not coprime to $(a+b)/c$ can happen for at most one prime $p$, and in this case $p$ divides $a+b$. Indeed, if $d=\gcd(u,(a+b)/c)$, then 
multiplying across equation (\ref{eq:****}) by $abu$ and reducing the resulting equation modulo $d$, we get $p\equiv 0\pmod d$. This is possible only if $d$ is prime and $p=d$, so $p$ divides $a+b$. Hence, 
$$
(\log x)^4\le m\le \frac{p+a+b}{abu}\le \frac{2(a+b)}{abu}\le \frac{4}{au}\le 4,
$$
a contradiction for large $x$. Thus, we may assume that $u$ is coprime to $(a+b)/c$. Then the number of such primes $p\le x$ is therefore
$$
\pi(x;abu,d^*)\ll \frac{x}{\varphi(abu)\log(x/(abu))}\ll \frac{x}{\varphi(abu)\log x},
$$
where the last inequality follows because $abu\le x^{1-\delta}$. Summing over $a,b,c$ and $u$, we get that the number of such situations is 
\begin{eqnarray*}
S_1 & := & \sum_{p\le x} f_1(p)\\
&  \ll & \frac{x}{\log x}\sum_{\substack{(a,b)=1 \\ ab\le x}} \sum_{c\mid a+b} \sum_{u\le x/ab} \frac{1}{\varphi(abu)}\\
& \ll & 
\frac{x}{\log x}\sum_{\substack{(a,b)=1\\  ab\le x}}\frac{\tau(a+b)}{\varphi(a)\varphi(b)}\sum_{u\le x/ab}\frac{1}{\varphi(u)}.
\end{eqnarray*}
The inner sum is $\le \sum_{u\le x} 1/\varphi(u)\ll \log x$. Thus,
$$
S_1 \ll x\sum_{\substack{(a,b)=1\\ ab\le x}} \frac{\tau(a+b)}{\varphi(a)\varphi(b)}.
$$
We use the fact that 
$$
\frac{1}{\varphi(n)}\ll \frac{\sigma(n)}{n^2}=\sum_{d\mid n} \frac{1}{dn}.
$$
With this, and writing $a=d_1u,~b=d_2v$ whenever $d_1,~d_2$ are divisors of $a$ and $b$ respectively, we get that the above quantity is 
\begin{eqnarray}
\label{eq:20}
S_1 \ll  x\sum_{\substack{(a,b)=1\\ ab\le x}} \frac{\tau(a+b)}{\varphi(a)\varphi(b)} & \ll & x\sum_{\substack{(a,b)=1\\ ab\le x}} \tau(a+b)\sum_{d_1\mid a}\sum_{d_2\mid b} \frac{1}{d_1ad_2b}\nonumber\\
& \ll & x \sum_{\substack{d_1\le x\\ d_2\le x\\ (d_1,d_2)=1}} \sum_{\substack{u\le x/d_1\\ (u,d_2)=1}} \sum_{\substack{v\le x/d_2\\ (v,d_1u)=1}} \frac{\tau(d_1u+d_2v)}{d_1^2d_2^2 uv}.
\end{eqnarray}
Proposition 7.6 in \cite{ElsTao} shows that uniformly in $A,~B,~C,~D$ all larger than $1$, we have
$$
\sum_{\substack{a\le A, b\le B, c\le C, d\le D\\ (ab,cd)=1}} \tau(ab+cd)\ll ABCD \log(A+B+C+D).
$$
Writing $A=2^i,~B=2^j,~C=2^k,~D=2^{\ell}$ for $i,j,k,\ell$  integers in $[0, \log x/\log 2]$, we have that
$$
\sum_{\substack{A\le d_1\le 2A, B\le d_2\le 2B, C\le u\le 2C, D\le v\le 2D\\ \gcd(d_1u,d_2v)=1}} \frac{\tau(d_1u+d_2v)}{d_1^2d_2^2 uv}\ll \frac{\log x}{AB}.
$$
Summing this up over all $i,j,k,\ell$ in $[0,\log x/\log 2]$ and putting $m:=i+j$, we get that 
\begin{eqnarray}
\label{eq:21}
\sum_{\substack{d_1\le x\\ d_2\le x\\ (d_1,d_2)=1}} \sum_{\substack{u\le x/d_1\\ (u,d_2)=1}} \sum_{\substack{v\le x/d_2\\ (v,d_1)=1}} \frac{\tau(d_1u+d_2v)}{d_1^2d_2^2 uv} & \ll &  
\log x \sum_{0\le i,j,k,\ell\le  \log x/\log 2} \frac{1}{2^{i+j}}\nonumber\\
& \ll &  (\log x)^3 \sum_{0\le m\le 2\log x/\log 2} \frac{m}{2^m}\nonumber\\
& \ll & (\log x)^3.
\end{eqnarray}
Inserting \eqref{eq:21} into \eqref{eq:20}, we thus get that 
$$
S_1\ll x(\log x)^3.
$$
This was under the assumption that $abu\le x^{1-\delta}$. So, from now on we assume that $abu>x^{1-\delta}$. 

\medskip

{\bf Case 2.} \textit{Assume $abm\le x^{1-\delta}$}.

\medskip

Let $f_2(p)$ be the number of such pairs $(m,p)$. To count 
$$
S_2=\sum_{p\le x} f_2(p),
$$
we let $a,b,c$ be fixed, then fix $m$ such that $abm\le x^{1-\delta}$ and we need to count the number of primes $p$ such that $(p+(a+b)/c)/(abm)=u$ is an integer. The number of such primes is
$$
\pi(x,abm,d^*)\ll \frac{x}{\varphi(abm)\log(x/(abm))}\ll \frac{x}{\varphi(abm)\log x}.
$$
The last inequality above holds since $abm\le x^{1-\delta}$. Here, similar to the previous case, we put $d^*$ for the class of $-(a+b)/c$ modulo $abm$. Again, $(a+b)/c$ is coprime to $m$, for if not, as in the analysis of the previous case, we get that $p\mid a+b$, so that
$$
1\le u\le \frac{p+a+b}{abm}\le \frac{2(a+b)}{abm}\le \frac{4}{m}\le \frac{4}{(\log x)^4},
$$
which is false for large $x$. Now an argument similar to the one from Case 1 (just swap the roles of $u$ and $m$) leads to
$$
S_2\ll x(\log x)^3.
$$

We next comment on the sizes of $a,b,c$ relative to each other. As we saw, we have $a\le b$. If $a=b$, then since $(a,b)=1$, we have that $a=b=1$ so $c\in \{1,2\}$. Thus, $m\mid p+1$ or $m\mid p+2$. 
Hence, the number of such situations is 
$$
\le \sum_{p\le x} (\tau(p+1)+\tau(p+2))\ll x.
$$
From now on, $a<b$. If $c=a+b$, then $m\mid p+1$. So, the number of such situations is $\sum_{p\le x}\tau(p+1)=O(x)$. We also assume that $c\le (a+b)/2<b$. Thus, $b>\max\{a,c\}$. 
We write $(a+b)/c=t$ so that $b=ct-a$ and 
$$
(ct-a)aum=p+t.
$$
Thus,
\begin{equation}
\label{eq:30}
t(acum-1)-a^2um=p.
\end{equation}
Clearly, 
$$
t(acum-1)=p+a^2um\le p+abum=2p+t\le 3x. 
$$

\medskip

{\bf Case 3.} \textit{Suppose that $t\ge x^{\delta}$.}

\medskip 

It follows that $acum\ll x^{1-\delta}$. We fix $a,c,u,m$ and count the number of primes $p\le x$ given by the form \eqref{eq:30}.
This is the same as counting the number of primes in some arithmetical progression of ratio $acum-1$ of first  term $a^2um$ coprime to $acum-1$. Note that $a^2um$ and $acum-1$ are coprime. By the Siegel-Walfisz theorem, the number of such primes is 
$$
\pi(x,acum-1,a^2um)\ll \frac{x}{\varphi(acum-1)\log(x/(acum-1))}\ll \frac{x}{\varphi(acum-1)\log x}.
$$
For the right--most inequality above, we used the fact that $acum\ll x^{1-\delta}$. 
The constant implied by the last Vinogradov symbol above, as well as most of the ones from the previous cases, depend on $\delta$ but at the end we will fix $\delta$ so all such constants are in fact absolute. 
So, we the contribution of this situation is 
\begin{equation}
\label{eq:40}
S_{3}\ll \frac{x}{\log x}\sum_{acum\le x} \frac{1}{\varphi(acum-1)}.
\end{equation}
We need to estimate the last sum. We now use the formula
$$
\frac{1}{\varphi(n)}\ll \sum_{d\mid n} \frac{1}{dn},
$$
but we truncate it $d<n^{1/5}$. Indeed, 
$$
\frac{1}{\varphi(n)}\ll \sum_{\substack{d\mid n\\ d\le n^{1/5}}} \frac{1}{dn}+O\left(\frac{\tau(n)}{n^{1+1/5}}\right).
$$
Since $\tau(n)=n^{o(1)}$, it follows that the last term on the right hand side is certainly $O(n^{-1-1/6})=o(1/n)$, so it can be absorbed into the left--hand side. With this, we get
$$
\sum_{acum\le x}\frac{1}{\varphi(acum-1)}\le \sum_{d\le (acum)^{1/5}} \sum_{\substack{acum\le x\\ d\mid acum-1}} \frac{1}{d(acum-1)}.
$$
Fix $d$ and $acum$ such that $d<(acum)^{1/5}$. Then $acum-1\equiv 0\pmod d$. There are various possibilities according to which one of the four numbers $a,c,u,m$ is larger. Say $c\ge \max\{a,u,m\}$. Then 
$c\ge (acum)^{1/4}\ge d^{5/4}$. Fix $d,a,u,m$. Then the congruence $c(aum)\equiv 1\pmod d$ puts $c$ into a progression $c^*\pmod d$, where $c^*\in [1,d-1]$. Let $c=c^*+dt$ for some $t$. Since $c^*<d\le (acum)^{1/5}\le c^{4/5}$, we get that $t\ge 1$. Then 
$$
acum-1=(aum)(c^*+dt)-1=d\left((aumt)+\frac{aum c^*-1}{d}\right)
$$
and the last fraction above is a positive integer. Thus, we get that
$$
\sum_{\substack{(a,c,u,m)\\ c\ge \max\{a,u,m\}}} \frac{1}{\varphi(acum-1)}\ll \sum_{d\le (acum)^{1/5}} \frac{1}{d^2 aumt}.
$$
We now sum over $d,a,u,m,t$, getting 
\begin{eqnarray*}
\sum_{\substack{acum\le x\\ c\ge \max\{a,u,m\}}}\frac{1}{\varphi(acum-1)} & \ll & \left(\sum_{d}\frac{1}{d^2}\right)\left(\sum_{a\le x} \frac{1}{a} \right)\left(\sum_{u\le x}\frac{1}{u}\right)\left(\sum_{t\le x} \frac{1}{t}\right)\left(\sum_{m\le x}\frac{1}{m}\right)
\\
& \ll & (\log x)^4.
\end{eqnarray*}
A similar situation happens if one of the other $3$ variables $a,u,m$ is $\max\{a,c,u,m\}$. Thus, 
$$
\sum_{acum\le x} \frac{1}{\varphi(acum-1)}\ll (\log x)^4,
$$
which inserted into \eqref{eq:40} gives $S_{3}\ll x(\log x)^3$. 

\medskip

{\bf Case 4.} \textit{The remaining case.}

\medskip

Here, we assume that $abu> x^{1-\delta},~abm> x^{1-\delta}, ~t< x^{\delta}$. It then follows that $\max\{m,u,t\}\ll x^{\delta}$. Since $abu>x^{1-\delta}$ and $a<b$, we get that $b>x^{(1-2\delta)/2}$. Since $b<ct\le cx^{\delta}$, 
we get that $c\ge x^{(1-4\delta)/2}$. Taking $\delta:=1/10$, we get that $c\gg x^{0.3}$. We return to equation \eqref{eq:30}, which we write as 
$$
c(amut)-(t+a^2um)=p.
$$
Observe that 
$$
c(amut)=p+t+a^2um\ll x+abum\ll x.
$$
Thus, $amut\ll x/c\ll x^{0.7}$. Also $t+a^2um$ is coprime to $amut$. Indeed, for if not, then the only possibility is that $p$ is a prime factor of the  number $\gcd(amut,t+a^2um)$ which must divide $t$. Thus, 
$p\le x^{\delta}$. Thus, the number of such pairs is at most 
$$
\sum_{p\le x^{\delta}} A_3(p)\le x^{\delta+1/2+o(1)}\le x^{2/3}\quad {\textrm{ for large}}\quad  x. 
$$
Fix $a,m,u,t$. We apply the Siegel-Walfisz theorem to get that for fixed $a,m,u,t$, the number of such primes is of order at most
$$
 \pi(x,amut,-(t+a^2um))\ll \frac{x}{\varphi(amut) \log(x/(amut))}\ll \frac{x}{(\log x)\varphi(amut)}.
$$
The last inequality follows again because $amut\ll x^{1-\delta}$. We now sum up over all $a,u,t,m$ getting 
$$
S_4\ll \frac{x}{\log x}\left(\sum_{a\le x} \frac{1}{\varphi(a)}\right)\left(\sum_{m\le x} \frac{1}{\varphi(m)}\right)\left(\sum_{u\le x} \frac{1}{\varphi(u)}\right)\left(\sum_{t\le x} \frac{1}{\varphi(t)}\right)\ll x(\log x)^3.
$$
This finishes the problem for the Type I solutions.

\subsection{Type II solutions}

Here, from (\ref{eq:1_1}), we have 
$$
m=\frac{1+p(a+b)/c}{abu}.
$$
Again, $(a,b)=1$ and $u$ is coprime to $(a+b)/c$ otherwise the above $m$ is not an integer. The cases when $abu\le x^{1-\delta}$ similar as in the case of solutions of Type I, since 
then the number of such pairs is 
\begin{eqnarray*}
 & \ll &  \sum_{\substack{abu\le x^{1-\delta}\\ (a,b)=1}}\sum_{c\mid a+b} \pi(x,abu,\rho^*)\\
 & \ll & \sum_{\substack{abu\le x^{1-\delta}\\ (a,b)=1}} \frac{x\tau(a+b)}{\phi(abu)\log(x/abu)} \\
& \ll & \frac{x}{\log x}\sum_{abu\le x}\frac{\tau(a+b)}{\phi(a)\phi(b)\phi(u)}\\
& \ll & 
x(\log x)^3,
\end{eqnarray*}
as in the analysis of the Type I situations. The same comment applies when $abm\le x^{1-\delta}$. From now on, we assume that $abu>x^{1-\delta}$ and $abm>x^{1-\delta}$. 
We write $t:=(a+b)/c$ and note that $b=ct-a$. Hence,
$$
1+pt=1+p(a+b)/c=abum=a(ct-a)um=(acum t)-a^2um,
$$
so 
$$
acum-\frac{(a^2um+1)}{t}=p.
$$
This signals $t$ as a divisor of $a^2um+1$. Further, since
$$
m=\frac{1+pt}{abu}=\frac{1+pt}{a(ct-a)u}>\frac{p}{acu},
$$
it follows that $acum>p$. Since $b\ge a$, we have that $ct\ge 2a$, so $ct-a\ge ct/2$, therefore
$$
m=\frac{1+pt}{a(ct-a)u}\le \frac{2tp}{a(ct/2)u}=\frac{4p}{acu},
$$
showing that $acum\le 4p$. In particular, $acum\asymp p$. Let us fix $a,u,m$ and $t\mid \tau(a^2um+1)$. Then 
\begin{equation}
\label{eq:acum}
acum+\frac{a^2um+1}{t}=p
\end{equation}
determines $p$ uniquely in terms of $c$. Since $c$ can be chosen in at most $4x/(aum)$ ways and $t$ in at most $a^2um+1$ ways, we get that 
the number of possibilities is 
$$
\le x\sum_{aum\le x} \frac{\tau(a^2um+1)}{aum}.
$$
It follows easily from Proposition 1.4 and Corollary 7.4 in \cite{ElsTao} that 
$$
\sum_{A\le a\le 2A}\sum_{U\le u\le 2U} \sum_{M\le m\le 2M} \tau(a^2um+1)\le AUM(\log x)^2,
$$
whenever  $A,~U,~M$ are positive integers in $[1,x]$. 
It then follows that 
$$
\sum_{A\le a\le 2A} \sum_{B\le b\le 2B} \sum_{M\le c\le 2M} \frac{\tau(a^2mu+1)}{amu}\ll (\log x)^2.
$$
Summing this up over all $(A,U,M)=(2^i,2^j,2^k)$ with $i,j,k$ integers in $[0,\log x/\log 2]$, we get an upper bound of $O(x(\log x)^5)$.

One may ask which bound is closer to the truth. We believe the lower bound is closer to the truth. Indeed, in the upper bound quite likely we have an extra $\log x$ factor for the sum $\tau(a^2 mu+1)$ over 
$a$'s (see the Remark 1.5 in \cite{ElsTao}). In addition, we should be able to save an extra  factor of $\log x$ in the upper bound by imposing that the expression in the left--hand side of \eqref{eq:acum} is prime (in our sum, we only summed over those $t$ such that $t\mid a^2um+1$ and did not use the extra condition that $acum+(a^2um+1)/t$ is prime). Because of these two extra conditions which we did not fully exploit, we conjecture that 
in fact the estimate
$$
\sum_{p\le x} A_3^*(p)\ll x(\log x)^3
$$
holds and leave this as an open problem. 

\section{Acknowledgements}

We thank the referee for a careful reading of the manuscript. This paper was written during visits of F.~L. to Universit\`a Roma Tre, in March, 2019 and Max Planck Institute for Mathematics in September, 2019. F.~L. thanks these Institutions for their hospitality and support. In addition, F.~L. was  
supported in part by grant CPRR160325161141 from the NRF of South Africa, by grant RTNUM18 from CoEMaSS at Wits, and by grant no. 17-02804S of the Czech Granting Agency.  F.~P. was supported in part by G.N.S.A.G.A. from I.N.D.A.M.

\end{document}